\documentclass[12pt]{article}

\usepackage{graphicx,tikz,color,url}
\usepackage{amsmath,amssymb,latexsym}
\usepackage{MnSymbol}

\newtheorem{theorem}{Theorem}[section]

\newtheorem{lemma}[theorem]{Lemma}
\newtheorem{proposition}[theorem]{Proposition}

\newtheorem{conjecture}[theorem]{Conjecture}
\newtheorem{corollary}[theorem]{Corollary}

\newtheorem{problem}[theorem]{Problem}

\newcommand{\proof}{\noindent{\bf Proof.\ }}
\newcommand{\qed}{\hfill $\square$ \bigskip}
\newcommand{\cp}{\,\square\,}

\newcommand{\name}{weakly claw-free}
\newcommand{\supp}{{\rm Supp}}

\renewcommand{\gg}{\gamma_{g}}
\newcommand{\ggz}{\gamma_{Zg}}
\newcommand{\ggt}{\gamma_{tg}}
\newcommand{\ggl}{\gamma_{Lg}}
\newcommand{\ggll}{\gamma_{LLg}}
\newcommand{\ggs}{\gamma_{g}'}
\newcommand{\ggzs}{\gamma_{Zg}'}
\newcommand{\ggts}{\gamma_{tg}'}
\newcommand{\ggls}{\gamma_{Lg}'}
\newcommand{\gglls}{\gamma_{LLg}'}

\newcommand{\vertex}{\node[vertex]}
\tikzstyle{vertex}=[circle, draw, inner sep=0pt, minimum size=6pt]

\begin{document}

\title{Z-domination game}

\author{Csilla Bujt\'as$^{a}$\thanks{Email: \texttt{csilla.bujtas@fmf.uni-lj.si}} 
\and Vesna Ir\v si\v c$^{a,b}$\thanks{Email: \texttt{vesna.irsic@fmf.uni-lj.si}}
\and Sandi Klav\v zar $^{a,b,c}$\thanks{Email: \texttt{sandi.klavzar@fmf.uni-lj.si}}
}
\maketitle

\begin{center}
$^a$ Faculty of Mathematics and Physics, University of Ljubljana, Slovenia\\
\medskip

$^b$ Institute of Mathematics, Physics and Mechanics, Ljubljana, Slovenia\\
\medskip

$^c$ Faculty of Natural Sciences and Mathematics, University of Maribor, Slovenia\\
\medskip
\end{center}

\begin{abstract}
The Z-domination game is a variant of the domination game in which each newly selected vertex $u$ in the game must have a not yet dominated neighbor, but after the move all vertices from the closed neighborhood of $u$ are declared to be dominated. The Z-domination game is the fastest among the five natural domination games. The corresponding game Z-domination number of a graph $G$ is denoted by $\gamma_{Zg}(G)$. It is proved that the game domination number and the game total domination number of a graph can be expressed as the game Z-domination number of appropriate lexicographic products. Graphs with a Z-insensitive property are introduced and it is proved that if $G$ is Z-insensitive, then $\gamma_{Zg}(G)$ is equal to the game domination number of $G$. Weakly claw-free graphs are defined and proved to be Z-insensitive. As a consequence, $\gamma_{Zg}(P_n)$ is determined, thus sharpening an earlier related approximate result. It is proved that if $\gamma_{Zg}(G)$ is an even number, then $\gamma_{Zg}(G)$ is strictly smaller than the game L-domination number. On the other hand, families of graphs are constructed for which all five game domination numbers coincide. Graphs $G$ with $\gamma_{Zg}(G) = \gamma(G)$ are also considered and computational results which compare the studied invariants in the class of trees on at most $16$ vertices reported. 
\end{abstract}


\noindent
{\bf AMS Subj.\ Class.\ (2010)}: 05C57, 05C69, 91A43

\section{Introduction}
\label{sec:intro}

Along with the well-established domination game (introduced in~\cite{bresar-2010}, see also~\cite{james-2019, kinnersley-2013, klavzar-2019, xu-2018b}) and the total domination game (introduced in~\cite{henning-2015}, see also~\cite{bujtas-2018, bujtas-2016, dorbec-2016, henning-2017}), the so-called Z-domination game, L-domination game, and LL-domination game were recently introduced in~\cite{bresar-2019+}. In this way a set-up of possible domination games became naturally rounded. Besides being a key stone of this classification, the Z-domination game is also motivated by the Grundy domination number~\cite{bgmrr-2014}. More precisely, in~\cite{BBG-2017} it was proved that the sum of the Z-Grundy domination number of a graph and its zero forcing number~\cite{AIM} is equal to the order of the graph. Now, selecting the next vertex in a Z-sequence and the next vertex in a Z-domination game is done using the same principle. Lin~\cite{lin-2019} further extended this idea by relating four variants of the zero forcing number to four variants of the Grundy domination number.

 For a vertex $v$ of a graph $G$, its open and closed neighborhoods are respectively denoted by $N(v)$ and $N[v]$. Each of the above listed games is played by Dominator and Staller who alternately select a vertex from $G$. If Dominator is the first to play we speak of a D-game, otherwise we have an S-game. In the $i^{\rm th}$ move, the choice of a vertex $v_i$ is legal if for the vertices $v_1,\ldots,v_{i-1}$ chosen so far, the following hold:
    \begin{itemize}
        \item $N[v_i] \setminus \bigcup_{j=1}^{i-1}N[v_j]\not=\emptyset$, in
        the {\em domination game};
        \item $N(v_i) \setminus \bigcup_{j=1}^{i-1}N(v_j)\not=\emptyset$, in
        the {\em total domination game};
        \item $N(v_i) \setminus \bigcup_{j=1}^{i-1}N[v_j]\not=\emptyset$, in {\em the Z-domination game};
        \item $N[v_i] \setminus \bigcup_{j=1}^{i-1}N(v_j)\not=\emptyset$ and $v_i\neq v_j$ for all $j <i$, in the {\em L-domination game}; and
        \item $N[v_i] \setminus \bigcup_{j=1}^{i-1}N(v_j)\not=\emptyset$, in
        {\em the LL-domination game}.
    \end{itemize}
Each of the games ends if there are no more legal moves available. We assume throughout this paper that the games are played on isolate-free graphs. Under this condition, a domination game and a Z-domination game on $G$ end with the $i^{\rm th}$ move if $\bigcup_{j=1}^{i}N[v_j]=V(G)$, while a total, L-, and LL-domination game end when $\bigcup_{j=1}^{i}N(v_j)=V(G)$. In each defined version of the game  Dominator wishes to finish it as soon as possible, while Staller wishes to delay the end. If a D-game is played and both players play optimally, the length of the game, i.e., the total number of moves played during the game, is, respectively, the {\em game domination number} $\gg(G)$, the {\em game total domination number} $\ggt(G)$, the {\em game Z-domination number} $\ggz(G)$, the {\em game L-domination number}  $\ggl(G)$,  and the {\em game LL-domination number} $\ggll(G)$ of $G$. For the S-game the corresponding invariants are $\ggs(G)$, $\ggts(G)$, $\ggzs(G)$, $\ggls(G)$, and $\gglls(G)$.
 
In this paper we are interested in the game Z-domination number and its interplay with other related (game) domination invariants. In the next section we present relations between the five game domination numbers and recall a couple of results needed later. We demonstrate in Section~\ref{sec:gg-and-ggt-with-ggz} that $\gg$ and $\ggt$ of a given graph can be expressed as $\ggz$ of an appropriate lexicographic product. Then, in Section~\ref{sec:ggz=gg}, we are interested in the relation between $\ggz$ and $\gg$. We introduce Z-insensitive graphs and prove that if $G$ is Z-insensitive, then $\gamma_{Zg}(G) = \gamma_g(G)$ as well as $\gamma_{Zg}'(G) = \gamma_g'(G)$. We also introduce weakly claw-free graphs (which properly include claw-free graphs) and prove that each weakly claw-free graph is Z-insensitive. As a direct consequence we deduce the exact values of $\gamma_{Zg}(P_n)$, thus sharpening~\cite[Theorem 5.3]{bresar-2019+}. Moreover, $\ggz$ is determined for powers of cycles. In Section~\ref{sec:ggz=ggl} we relate $\ggz$ with $\ggl$.  We first list several families of graphs for which the two invariants are the same. In particular, combining results from this paper and an earlier paper we demonstrate that if $n \geq 2 m - 1$ and $m \geq 2$, then $\ggz(K_m \cp K_n) = \ggl(K_m \cp K_n) = \ggll(K_m \cp K_n) = 2m - 1$. On the other hand we prove that if $\ggz(G)$ is an even number, then $\ggz(G) <  \ggl(G)$. We also conjecture that $\ggz(T) <  \ggl(T)$ holds for an arbitrary tree of order at least $2$. In a brief Section~\ref{sec:ggz=gamma} we associate to an arbitrary graph $G$ on at least three vertices a graph $\widehat{G}$ such that $\ggz(\widehat{G}) = \gamma(\widehat{G}) = \frac{\gg(\widehat{G}) + 1}{2}$ holds. In the final section we report on our systematic computation of all the studied invariants in the class of trees on at most $16$ vertices. 

\section{Preliminaries}
\label{sec:preliminaries}

 Let $G$ be a graph that has minimum degree $\delta(G) \ge 1$. The order of $G$ will be denoted with $n(G)$. A set $S \subseteq V(G)$ is a \emph{dominating set} of $G$ if $N[S] = V(G)$ and $S$ is a \emph{total dominating set} of $G$ if $N(S) = V(G)$.   The minimum cardinality of a dominating set (resp., total dominating set) is the \emph{domination number} $\gamma(G)$ of $G$ (resp., \emph{total domination number} $\gamma_t(G)$). By the definition of the games, $\gamma(G)$ is a lower bound for $\ggz(G)$ and $\gg(G)$, while $\gamma_t(G)$ is a lower bound for each of $\ggt(G)$, $\ggl(G)$ and $\ggll(G)$.

The following result describes the basic relations between the game domination numbers. 

\begin{theorem} {\rm \cite[Theorem 3.1]{bresar-2019+}}
\label{thm:hierarchy}
If $G$ is a graph without isolated vertices, then
$$\ggz(G)\le \gg(G) \le \ggl(G)\le \ggll(G) \quad \text{and} \quad \ggz(G)\le \ggt(G) \le \ggl(G)\le \ggll(G)\,.$$
\end{theorem}

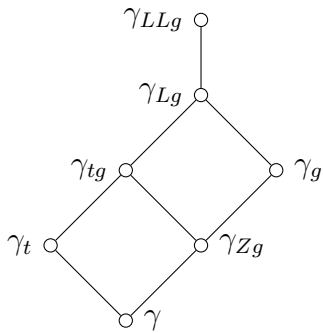
\begin{figure}[ht!]
    \begin{center}
        \begin{tikzpicture}[]
        \tikzstyle{vertex}=[circle, draw, inner sep=0pt, minimum size=5pt]
        \tikzset{vertexStyle/.append style={rectangle}}
        \vertex (1) at (0,0) [label=right:$\ggz$] {};
        \vertex (2) at (-1,1) [label=left:$\ggt$] {};
        \vertex (3) at (1,1) [label=right:$\gg$] {};
        \vertex (4) at (0,2) [label=left:$\ggl$] {};
        \vertex (5) at (0,3) [label=left:$\ggll$] {};
        \vertex (6) at (-1,-1) [label=right:$\gamma$] {};
        \vertex (7) at (-2,0) [label=left:$\gamma_t$] {};
        \path
        (1) edge (2)
        (1) edge (3)
        (2) edge (4)
        (3) edge (4)
        (4) edge (5)
        (1) edge (6)
        (6) edge (7)
        (7) edge (2);
        \end{tikzpicture}
    \end{center}
    \caption{Relations between the five versions of the game domination number, the domination number, and the total domination number.} 
\label{fig:Hasse}
\end{figure}

Note that $\gamma_t(P_6) = 4$ and $\gg(P_6) = \ggz(P_6) = 3$ (cf.\ Corollary~\ref{thm:paths}), and that $\gamma_t(P_2\cp P_3) = 2$ and $\gg(P_2\cp P_3) = \ggz(P_2\cp P_3) = 3$. Hence $\gamma_t$ is incomparable with both $\gg$ and $\ggz$. Note also that by definition, $\gamma(G)\le \ggz(G)$. Moreover, as $\gg$ and $\ggt$ are incomparable~\cite{henning-2015}, the Hasse diagram on the above invariants with respect to the relation $\le$ is as shown in Fig.~\ref{fig:Hasse}. 

Let $\supp(G)$ denote the set of support vertices of $G$, that is, vertices which have at least one neighbor of degree $1$. If $\supp(G)$ forms a dominating set of $G$, then it is called a \textit{supportive dominating set}. Note that if $\supp(G)$ is a supportive dominating set of a graph $G$, then $|\supp(G)| = \gamma(G)$. 

\begin{theorem}\label{thm:enough-pendant} {\rm \cite[Theorem 3.1]{xu-2018}} 
Let $G$ be a connected graph of order at least $3$. If $G$ has a supportive dominating set and there are at least $\lceil\log_2 \gamma(G)\rceil+1$ pendant vertices adjacent to each vertex of $\supp(G)$, then  $\gg(G) = 2\gamma(G) - 1$. 
\end{theorem}

Two vertices, $u$ and $v$, are \emph{(true) twins} in $G$, if $N_G[u]=N_G[v]$, and they are \emph{false twins} if $N_G(u)=N_G(v)$.

\begin{lemma}\label{lem:twin} {\rm \cite[Proposition 1.4]{bresar-2017}}
	If $G$ is a graph and $u, v \in V(G)$ are twins, then $\gg(G) = \gg(G-v)$.
\end{lemma}

The following lemma can be easily derived from the fact that if $u$ and $v$ are false twins, then at most one of them can be played during a total domination game. Alternatively, it follows from \cite[Lemma 12]{henning-2018}.

\begin{lemma}\label{lem:falsetwin}
	If $G$ is a graph with $\delta(G) \ge 1$ and $u, v \in V(G)$ are false twins, then $\ggt(G) = \ggt(G-v)$.
\end{lemma}

If $G$ is a graph and $A\subseteq V(G)$, then $G|A$ denotes a {\em partially dominated graph} meaning that when a game is played on $G|A$, the vertices from $A$ need not be dominated but they are allowed to be played provided they are legal moves.

\section{$\gg$ and $\ggt$ expressed with $\ggz$}
\label{sec:gg-and-ggt-with-ggz}

In this section we prove that the game domination number and the total game domination number of a graph $G$ can be expressed with the game Z-domination number of an appropriate lexicographic product. For this sake recall that the \emph{lexicographic product} $G \circ H$ of graphs $G$ and $H$ is a graph with vertices $V(G) \times V(H)$, where vertices $(g_1, h_1)$ and $(g_2, h_2)$ are adjacent if $g_1g_2 \in E(G)$, or $g_1=g_2$ and $h_1h_2 \in E(H)$. We first give the connection for the total game domination number, where $\overline{X}$ denotes the complement of a graph $X$. 

\begin{theorem}
	\label{thm:lexi_complement_Zg}
	If $n \geq 2$ and $\delta(G) \ge 1$, then $\ggt(G) = \ggz(G \circ \overline{K}_n)$.
\end{theorem}

\proof
Let $V(G)=\{v_1, \dots ,v_\ell \}$ and $V(\overline{K}_n)=[n]$. We will use the notation  $F=G \circ \overline{K}_n$ and $V_i=\{(v_i,j) \mid j\in [n]\}$ for $i\in [\ell]$. Since each $V_i$ contains false twin vertices, a repeated application of Lemma~\ref{lem:falsetwin} implies $\ggt(F)=\ggt(G)$. By Theorem~\ref{thm:hierarchy}, we have $\ggz(F) \le \ggt(F)$. Thus it remains to prove $\ggz(F) \ge \ggt(G)$. 

Consider a Z-domination game on $F$ and observe that from each $V_i$ at most one vertex can be played during the game.  Indeed, after playing  a vertex $(v_i, j)$, all vertices from the open neighborhood $N_F((v_i, j'))$ become dominated for every $j' \in [n]$.

To prove $\ggz(F) \ge \ggt(G)$, we define two parallel games. The real game is a Z-domination game on $F$ where Dominator plays optimally; the other game is a total domination game on $G$ which is imagined by Staller. In the latter one Staller plays optimally. If Dominator plays a vertex $(v_i,j)$ in the real game, Staller interprets this move as $v_i$ in the imagined game. Then, Staller replies in the imagined total domination game with an optimal move $v_s$ and copies it into the real game as $(v_s,1)$. We will show that the moves remain legal when copied to the other game according to the described rules.

Let $D_R$ be the set of vertices that have been played in the real game until a point. Define 
$$A_F(D_R)=\{v_i\in V(G) \mid V_i \subseteq N_F[D_R]\}.$$
That is, $v_i \in A_F(D_R)$ if and only if every vertex from $V_i$ has been dominated in the real game. We prove that the following is true after every move.

\medskip\noindent
{\bf Property A}: After a move and its interpretation in the other game, $ A_F(D_R)= N_G(D_I)$ holds for $D_R$ and $D_I$ which are the sets of played vertices in the real and imagined games respectively.

\medskip  
The equality in Property A clearly holds at the beginning of the game. Suppose that $ A_F(D_R)= N_G(D_I)$ holds before a move of Dominator. Then, if he selects a vertex $(v_i,j)$ in the real Z-game, it has to dominate at least one vertex $(v_s,j')$ with $s\neq i$. It follows that $v_iv_s \in E(G)$ and $v_s \notin A_F(D_R)$. By  $ A_F(D_R)= N_G(D_I)$, we conclude that $v_s$ is not totally dominated in the imagined game and hence, the interpretation $v_i$ is a legal move in the imagined game. Recall that as $n \geq 2$, the vertex $(v_i,j)$ must be the first (and only) vertex played from $V_i$ and therefore, $v_i$ is not added to $A_F(D_R)$ when $(v_i,j)$ is played. In general, a vertex $v_s$ is added to $ A_F(D_R)$ after this move if and only if $v_iv_s \in E(G)$ and $v_s \notin  A_F(D_R)$. By the hypothesis $ A_F(D_R)= N_G(D_I)$, this exactly means that $v_s$ becomes totally dominated in the imagined game and added to $N_G(D_I)$. This proves that Dominator's move maintains Property A.

Consider next a move $v_i$ of Staller in the imagined game. With this move she totally dominates a vertex $v_p$ if and only if $v_iv_p \in E(G)$ and $v_p \notin N_G(D_I)$ for the set $D_I$ of vertices which have been played so far. By the hypothesis $ A_F(D_R)= N_G(D_I)$, not all vertices from $V_p$ are dominated in the real game before this move and, therefore, playing $(v_i, 1)$ is a legal move in the real game. Observe again that $(v_i,1)$ is the first vertex which is played from $V_i$. By the same reasoning as before, it can be proved that $v_p$ is added to $N_G(D_I)$ after the move $v_i$ if and only if it is added to $ A_F(D_R)$ after the move $(v_i, 1)$. This proves that Staller's move maintains Property A.

As Property A also holds at the end of the game and the real game finishes when  $ A_F(D_R)= V(G)$ is achieved and the imagined game finishes when  $N_G(D_I)=V(G)$ holds, the two games end in the same number of moves, say $t$.  Since Dominator is playing optimally in the real game and Staller in the imagined game, we have $ \ggz(F) \ge t \ge \ggt(G)$. This completes the proof.
\qed

For the game domination number we have the following parallel result, the proof of which will be given after Theorem~\ref{thm:Z-config}. 

\begin{proposition}
	\label{prp:lexi_complete_Zg}
	If $n \geq 2$, then $\gg(G) = \ggz(G \circ K_n)$.
\end{proposition}

\section{On the equality $\ggz = \gg$}
\label{sec:ggz=gg}

In this section we give a sufficient condition on a graph $G$ such that $\ggz(G) = \gg(G)$ holds. We then define a class of weakly claw-free graphs which contains claw-free graphs as a proper subclass. It is proved that weakly claw-free graphs have the property $\ggz(G) = \gg(G)$, and as a corollary we determine the game Z-domination number of paths and of powers of cycles. The first of these two consequences sharpens~\cite[Theorem 5.3]{bresar-2019+} where the game Z-domination number of paths was determined up to a constant $c$ with $|c|\le 2$. 
 
We say that a partially dominated graph $G|A$ has a \emph{Z-configuration} if there exists an undominated vertex $v$ such that all its neighbors are dominated and for every neighbor $u \in N(v)$ there exists an undominated neighbor different from $v$. More formally, a Z-configuration exists in $G|A$  if there is a vertex $v\in V(G)\setminus A$ such that $N(v) \subseteq A$  and we have $|N(u) \setminus A| \ge 2$ for each $u \in N(v)$.  We say that an isolate-free graph $G$ is \emph{Z-insensitive}  if the partially dominated graph $G|N[D]$ has no Z-configuration for all $D \subseteq V(G)$.

\begin{theorem}
	\label{thm:Z-config}
	If $G$ is Z-insensitive, then $\ggz(G) = \gg(G)$ and $\ggz'(G) = \gg'(G)$.
\end{theorem}

\proof
By Theorem~\ref{thm:hierarchy}, we only need to prove that $\ggz(G) \ge \gg(G)$. For this sake consider the following two games on $G$.

The real game is the Z-domination game played on $G$. In this game Dominator is playing optimally but Staller maybe not. At the same time Staller is imagining that the usual domination game is also played on $G$, in which Staller is playing optimally, but Dominator maybe not. Each move of Dominator from the real game is copied by Staller into the imagined game. Then, in the latter game, Staller replies optimally (with respect to the usual domination game that is played in the imagined game) and copies her move into the real game, provided the move is legal there. Otherwise Staller plays some other vertex, the selection of which is described below. We prove that the following property holds after every move of Staller. 

\medskip\noindent
{\bf Property B}: after Staller plays a move in the imagined game and a related (not necessarily the same) move in the real game, the set of vertices dominated is the same in both games. 

\medskip
Note that if Property B holds after a move of Staller, then the next move of Dominator played in the real game is a legal move in the imagined game, so that Staller can indeed copy the moves of Dominator into the imagined game. 

If Staller can copy her move into the real game (that is, if the move of Staller played in the imagined game is a legal move in the real Z-game), then Property B is clearly maintained. Suppose now that at some stage of the imagined game Staller selects a vertex $u$ which is not a legal move in the real, Z-domination game. This can only happen if in the real game $u$ is not yet dominated at this stage of the game but all its neighbors are. Let  $D$ denote the set of vertices which have been played in the real game so far. Since $G|N[D] $ does not have a Z-configuration, there exists a vertex $w \in N(u)$ such that $|N(w) \setminus N[D]| < 2$ and hence, as $w \in N[D]$ and $u \notin N[D]$, we have $N[w]\setminus N[D]= \{u\}$. Then the move  $w$  is clearly legal in the real (Z-domination) game and Staller can play it instead of $u$. The relation $N[D\cup \{w\}]=N[D\cup \{u\}]= N[D]\cup\{u\}$ ensures that Staller maintains Property B with this selection.

We have thus proved that Property B holds after each move. As a consequence, the real game and the imagined game ends in the same number of moves, say $t$. Since in the real, Z-domination game, Dominator is playing optimally, we have $\ggz(G)\ge t$. On the other hand, in the imagined game Staller is playing optimally  and thus, $\gg(G)\le t$. We conclude that $\ggz(G)\ge t\ge \gg(G)$.

The equality $\ggz'(G) = \gg'(G)$ is proved with the same reasoning.
\qed

As a first application of Theorem~\ref{thm:Z-config}, we give a short proof of Proposition~\ref{prp:lexi_complete_Zg}.

\medskip \noindent
\textbf{Proof of Proposition~\ref{prp:lexi_complete_Zg}.}
Note that vertices in each $K_n$-layer are twins. Thus by Lemma~\ref{lem:twin} it holds $\gg(G \circ K_n) = \gg(G)$. As $n \geq 2$, every vertex has a twin, thus for no vertex it can hold that it is undominated and all its neighbors are dominated. Hence $G \circ K_n$ is Z-insensitive and the equality $\gg(G \circ K_n) = \ggz(G \circ K_n)$ follows by Theorem~\ref{thm:Z-config}. 
\qed

We say that a graph $G$ is {\em \name} if every vertex $u\in V(G)$ has a neighbor that is not the center of a claw. Clearly, the class of \name\ graphs properly contains the class of claw-free graphs. Equally obvious is that if every vertex of a graph $G$ has a neighbor of degree at most $2$, then $G$ is \name. Observe that \name\ graphs have no isolated vertices.

\begin{theorem}
\label{thm:locally-claw-free}
If $G$ is a \name\ graph, then $\ggz(G) = \gg(G)$ and $\ggz'(G) = \gg'(G)$. 
\end{theorem}

\proof
We prove that every \name\ graph is Z-insensitive. Consider a set $D \subseteq V(G)$ and the partially dominated graph $G|N[D]$. Suppose for a contradiction that we have a Z-configuration at a vertex $u$. Then, $u$ is the only vertex from $N[u]$ which is not dominated by $D$. Further, every neighbor of $u$ has at least one neighbor different from $u$ that is not dominated by $D$. Since $G$ is \name, among the neighbors of $u$ there exists a vertex $w$ that it is not the center of a claw. Let $w'$ be a neighbor of $w$ different from $u$ that is not  dominated. Note that $u$ is not adjacent to $w'$ since by our assumption all the neighbors of $u$ are dominated. Since $w$ is  dominated but $u$ is not, there is a neighbor $w'' \in D$ which dominates $w$. Clearly,  $w''\ne w'$ and $w''\ne u$. Now, $w''$ is not adjacent to $u$ (for otherwise $u$ would be dominated), hence $w''$ must be adjacent to $w'$, for otherwise $w$ would be the center of a claw induced by $w$, $u$, $w'$, and $w''$. But $w''$ being adjacent to $w'$ means that $w'$ is dominated, a contradiction. Hence, $G$ is Z-insensitive and Theorem~\ref{thm:Z-config} implies the statement.
\qed
 
In~\cite[Theorem 5.3]{bresar-2019+} it was proved that for every positive integer $n$ there exists a constant $c_n$ such that $\ggz(P_n) = \frac{n}{2} + c_n$ holds with $|c_n|\le 2$. Combining Theorem~\ref{thm:locally-claw-free} with the known formula for $\gg(P_n)$, see~\cite[Theorem 2.4]{kosmrlj-2017}, we can strengthen~\cite[Theorem 5.3]{bresar-2019+} as follows.   

\begin{corollary}
\label{thm:paths}
If $n\ge 2$, then
$$\ggz(P_n) = 
\left\{
\begin{array}{ll}
\left\lceil\frac{n}{2}\right\rceil-1; & n\bmod (4) = 3\,, \\ \\
\left\lceil\frac{n}{2}\right\rceil; & otherwise\,.
\end{array}
\right.$$
\end{corollary}

Similarly, we get the exact value for $\ggz'(P_n)$. Since the powers of cycles are also claw-free graphs, Theorem~\ref{thm:Z-config} and the results~\cite[Theorem 9]{bujtas-2015} for $\gg(C_N^n)$ and $\gg'(C_N^n)$ can be rewritten for $\ggz(C_N^n)$ and  $\ggz'(C_N^n)$, respectively. Here, we only state the exact result for the D-game. 

\begin{corollary} \label{thm:power}
	For every $n\geq 1$ and $N \ge 3$,
	\begin{equation*}
	\ggz(C_N^n)=  \left\{
	\begin{array}{lll}
	\left\lceil \frac{N}{n+1}\right\rceil; & N\bmod (2n+2) \in \{0,1,\ldots, n+1\}\,,\\\\
	\left\lceil \frac{N}{n+1}\right\rceil-1; & N\bmod (2n+2)\in \{n+2,\ldots, 2n+1\}\,.
	\end{array} \right.
	\end{equation*}
\end{corollary}

Since $\gamma(G) \le \ggz(G) \le \gg(G)$, additional examples of graphs $G$ for which $\ggz(G) = \gg(G)$ holds are the graphs $G$ for which we have $\gamma(G) = \gamma_g(G)$. Trees with this property have been characterized in~\cite{nadjafi-2016}.

\section{On the equality $\ggz = \ggl$}
\label{sec:ggz=ggl}

In Section~\ref{sec:ggz=gg} we have found large classes of graphs $G$ with $\ggz(G) = \gg(G)$. Considering the Hasse diagram in Fig.~\ref{fig:Hasse}, we are next interested in graphs $G$ with $\ggz(G) = \ggl(G)$ or even with $\ggz(G) = \ggll(G)$, which arise from~\cite[Problem 6.4]{bresar-2019+} and~\cite[Problem 6.2]{bresar-2019+}, respectively. The only example found in~\cite{bresar-2019+} with equal Z-domination and LL-domination game numbers is the family of Cartesian product graphs $G \cp K_{1,k}$, where $G$ is a connected graph with $n(G) \geq 2$ and $k \geq 2 n(G)$~\cite[Proposition 6.1]{bresar-2019+}. (Recall that the {\em Cartesian product} $G\cp H$ of graphs $G$ and $H$ has the vertex set $V(G)\times V(H)$, vertices $(g,h)$ and $(g',h')$ being adjacent if either $gg'\in E(G)$ and $h=h'$, or $g=g'$ and $hh'\in E(H)$.) In this case it holds $\ggz(G \cp K_{1,k}) = \ggll(G \cp K_{1,k}) = 2 n(G) - 1$. We next present another infinite family of graphs for which all five game domination numbers are the same. As in the already known example, we will apply the Cartesian product of graphs. 

\begin{proposition}
	\label{prp:cartesian_product}
	If $n \geq 2 m - 1$ and $m \geq 2$, then $$\ggz(K_m \cp K_n) = \ggl(K_m \cp K_n)= \ggll(K_m \cp K_n) = 2m - 1\,.$$
\end{proposition}

\proof
As $n \geq 2 m - 1$, it follows by~\cite[Proposition 5.1]{bresar-2017} that $\gg(K_m \cp K_n) = 2m - 1$. As the graph $K_m \cp K_n$ is claw-free, we have $\gg(K_m \cp K_n) = \ggz(K_m \cp K_n)$ by Theorem~\ref{thm:locally-claw-free}.

It follows from Theorem~\ref{thm:hierarchy} that $\gg(K_m \cp K_n) \leq \ggll(K_m \cp K_n)$. On the other hand, we have $\ggll(K_m \cp K_n) \leq 2 \gamma_t(K_m \cp K_n) - 1 \leq 2m-1 = \gg(K_m \cp K_n)$ by~\cite[Proposition 4.1.(iii)]{bresar-2019+}.
\qed

In view of Fig.~\ref{fig:Hasse}, Proposition~\ref{prp:cartesian_product} implies that also $\gg(K_m \cp K_n) = \ggt(K_m \cp K_n) = \ggz(K_m \cp K_n) = 2m - 1$ holds whenever $n \geq 2 m - 1$ and $m \geq 2$. Hence, for these Hamming graphs all five game domination numbers coincide. 

There are additional families with this property. Let $G_{m,n}$, $m, n \geq 3$, be a graph obtained from disjoint copies of $K_m$ with vertex set $\{u_1, \ldots, u_m\}$ and $K_n$ with vertices $v_1, \ldots, v_n$, by adding the edges $u_1 v_1$ and $u_2 v_2$. It can easily be seen that $\ggz(G_{m,n}) \geq 3$, and on the other hand, we have $\ggll(G_{m,n}) \leq 2 \gamma_t(G_{m,n}) - 1 = 3$. Thus $\ggz(G_{m,n}) = \ggl(G_{m,n})= \ggll(G_{m,n}) = 3$.

Note that in all the examples of graphs $G$ given above for which $\ggz(G) = \ggl(G)$ holds, $\ggz(G)$ is an odd number. This is not a coincidence as the next result asserts. 

\begin{theorem}
	\label{thm:even_ggz}
If $\ggz(G)$ is an even number, then $\ggz(G)+1 \le \ggl(G)$.
\end{theorem}
\proof
Consider the following two parallel games on $G$. The real game is an L-domination game on $G$ in which Dominator plays with an optimal strategy. The imagined game is a Z-domination game on $G$ where Staller plays optimally. If Dominator plays a vertex in the real L-game, it is copied or interpreted (complying with some rules) in the imagined Z-game. Then Staller plays an optimal response which is copied or interpreted again in the real game. Let $D_R$ and $D_I$ denote the set of vertices which have been played until a point  in  the real and in the imagined game, respectively. We prove that the following property can be maintained after every move of the game.

\medskip\noindent
{\bf Property C}: After a move and its interpretation in the other game, $N[D_R] \subseteq N[D_I]$ holds.

\medskip  
Property C clearly holds with $D_R=D_I= \emptyset$ at the beginning of the game. Suppose first that $N[D_R] \subseteq N[D_I]$ holds before  Dominator plays a vertex $v$ in the L-game. If $v$ is also a legal move in the imagined Z-game, then it is copied there and as $N[D_R \cup \{v\}] \subseteq N[D_I \cup\{v\}]$ holds, Property C is maintained. If $v$ is not a legal move in the Z-game, then we have two cases. If $v$ and all neighbors of $v$ are already contained in $N[D_I]$, then $N[D_R \cup \{v\}] \subseteq N[D_I]$  and thus any legal move $v'$ in the imagined game can be the interpretation of the move $v$, hence Property C remains valid. In the other case, $N[v]\setminus N[D_I]=\{v\}$ and an arbitrary neighbor $v'$ of $v$ in the imagined game is a legal move which maintains Property C.

Consider next a move $u$ of Staller in the imagined game. Since it is a legal move in the Z-game, there exists a neighbor $u'$ of $u$ such that $u' \notin N[D_I]$. If Property C holds before this move, it implies $u' \notin N(D_R)$ and $u$ is a legal move in the real game. Then, if it was not the last move in the imagined game, we copy the move $u$ into the real game that maintains Property C. If the move $u$ finishes the imagined game, we may interpret it as a move $u'$ in the L-game. It is a legal move there indeed, as $u' \notin N[D_R]$ implies that the set $N[u']\setminus N(D_R)$ is not empty and $u' \notin D_R$. Then, as $u'$ is not totally dominated by $D_R \cup \{u'\}$, the real game does not end with this move.

Let $t_R$ and $t_I$ denote the number of moves in the real and imagined game, respectively. Since Dominator plays according to an optimal strategy in the real game,  $t_R \le \ggl(G)$ holds. Similarly, we have $t_I \ge \ggz(G)$. By Property C, we have $N(D_R) \subseteq N[D_I]$ that ensures $t_I \le t_R$. If the imagined game finishes with a move of Dominator, then $t_I$ is odd and, since $\ggz(G)$ is even by assumption, we have $\ggz(G) <t_I \le t_R \le \ggl(G)$ which establishes the statement of the theorem. In the other case the imagined game finishes with a move of Staller and, as we have seen, the real game is strictly longer than the imagined one. This gives $t_I <t_R$ and then, $\ggz(G)\le t_I < t_R \le \ggl(G)$  proves the statement.
\qed

We have seen that there are many graphs $G$ for which $\ggz(G) = \ggl(G)$ holds, and that in such cases $\ggz(G)$ must be odd. We conjecture that there are no such examples in the class of trees. 

\begin{conjecture}
	\label{conj:ZL-trees}
	If $T$ is a tree with $n(T) \geq 2$, then $\ggz(T) < \ggl(T)$.
\end{conjecture}

Since $\ggl(G)\le \ggll(G)$ holds for every graph $G$, Conjecture~\ref{conj:ZL-trees} strengthens \cite[Conjecture 6.3]{bresar-2019+} which asserts that $\ggz(T) < \ggll(T)$ holds for every tree $T$. We have checked by computer that Conjecture~\ref{conj:ZL-trees} holds for all trees $T$ with $2\le n(T)\le 18$.

\section{On the equality $\ggz = \gamma$}
\label{sec:ggz=gamma}

\medskip
From Theorem~\ref{thm:hierarchy} we in particular know that $\ggz(G) \le \gg(G)$ holds for every graph $G$. On the other hand, $\ggz(G)$ can also be bounded from below in view of $\gg(G)$ as follows. From~\cite[Theorem 1]{bresar-2010} we know that $\gamma(G)\le \gg(G) \le 2 \gamma(G)-1$, while from~\cite[Proposition 4.1(i)]{bresar-2019+} we also have $\gamma(G) \le \ggz(G)$. Hence: 
 
\begin{equation*}
\label{eq:chain}
\frac{\gg(G)+1}{2} \le \gamma(G) \le \ggz(G)\,.
\end{equation*}

For an arbitrary connected graph we can construct a related graph for which both equalities hold above. Let $G$ be a connected graph of order at least $3$, and let $\widehat{G}$ be the graph obtained from $G$ as follows. Add a vertex $w$ and connect it with an edge to every vertex of $G$. Then for each vertex $u\in V(G) \cup \{w\}$ add $\lceil\log_2 (n(G) + 1)\rceil+1$ pendant vertices adjacent to $u$. 

\begin{proposition}
	\label{prp:gamma_Zg}
	If $G$ is a graph with $n(G) \geq 3$, then $\ggz(\widehat{G}) = \gamma(\widehat{G}) = \frac{\gg(\widehat{G}) + 1}{2}$.
\end{proposition}

\proof
We use the notation introduced in the definition of the graph $\widehat{G}$. Note that $V(G) \cup \{w\}$ is a supportive dominating set, hence by Theorem~\ref{thm:enough-pendant}, $\gg(\widehat{G}) = 2\gamma(\widehat{G}) - 1$ and so $\gg(\widehat{G}) = 2n(G) + 1$. Let now Z-game be played on $\widehat{G}$ and let Dominator play $w$ in his first move. This move forces both players to play the remaining vertices from $\supp(\widehat{G})$ (that is, the vertices from $V(G)$) in the rest of the game. Hence $\ggz(\widehat{G}) \le n(G) + 1$. Since  $n(G) + 1 = \gamma(\widehat{G}) \le \ggz(\widehat{G})$ we have
$$\ggz(\widehat{G}) = n(G) + 1 = \frac{(2n(G)+1)+1}{2} = \frac{\gg(\widehat{G})+1}{2}\,,$$
and we are done. 
\qed

\section{Computations on trees}
\label{sec:conlcude}

In this section we present computational results in which $\ggz$ is compared to other invariants from Fig.~\ref{fig:Hasse} on the class of trees. 

The results are collected in Table~\ref{tbl:trees}. In the second column the number $\# T$ of trees $T$ of order $n$ is listed. In the next three columns we compare $\ggz$ with the three game domination numbers that are adjacent to $\ggz$ in Fig.~\ref{fig:Hasse}. More precisely, the number of trees has been computed for which $\ggz$ equals to one of these three invariants, respectively. We have already compared $\ggz$ with $\ggl$ at the end of Section~\ref{sec:ggz=ggl} and found out that such a column would contain only zeros. The remaining comparison is between $\ggz$ and $\gamma_t$. Since these two invariants are in general incomparable, we present in the last two columns that data for $\ggz > \gamma_t$ and $\ggz < \gamma_t$. From these data it follows that $\ggz$ and $\ggt$ are incomparable already on the class of trees. 

\begin{table}[ht!]
    \begin{center}
        \begin{tabular}{ | c || c | c | c | c | c | c |}
            \hline
            $n$ & $\!\!\#T\!\!$ & $\!\!\#T: \ggz = \gg\!\!$ & $\!\!\#T: \ggz = \ggt\!\!$ & $\!\!\#T: \ggz = \gamma\!\!$ & $\!\!\#T: \ggz > \gamma_t\!\!$ & $\!\!\#T: \ggz < \gamma_t\!\!$ \\ \hline \hline
            4 & $2$ & $2$ & $0$ & $2$ & $0$ & $1$ \\ \hline
            5 & $3$ & $3$ & $1$ & $2$ & $0$ & $1$ \\ \hline
            6 & $6$ & $5$ & $1$ & $4$ & $0$ & $2$ \\ \hline
            7 & $11$ & $10$ & $3$ & $6$ & $0$ & $3$ \\ \hline
            8 & $23$ & $19$ & $3$ & $11$ & $0$ & $6$ \\ \hline
            9 & $47$ & $40$ & $7$ & $16$ & $1$ & $8$ \\ \hline
            10 & $106$ & $84$ & $11$ & $29$ & $5$ & $21$ \\ \hline
            11 & $235$ & $186$ & $21$ & $47$ & $20$ & $41$ \\ \hline
            12 & $551$ & $412$ & $38$ & $84$ & $60$ & $103$ \\ \hline
            13 & $1301$ & $974$ & $75$ & $137$ & $189$ & $224$ \\ \hline
            14 & $3159$ & $2277$ & $141$ & $237$ & $559$ & $563$ \\ \hline
            15 & $7741$ & $5456$ & $277$ & $387$ & $1624$ & $1328$ \\ \hline
            16 & $19320$ & $13095$ & $539$ & $647$ & $4571$ & $3336$ \\ \hline
        \end{tabular}

    \end{center}

    \caption{The Z-game domination number on trees compared to other game domination numbers. The number of all trees $T$ of order $n$ is denoted by $\# T$.}
    \label{tbl:trees}
\end{table} 

Each of the five columns of Table~\ref{tbl:trees} leads to the question how the column continues. We explicitly state a corresponding problem for the first column as follows. 

\begin{problem}
Investigate the asymptotic behavior of the number of trees $T$ for which $\ggz(T) = \gg(T)$ holds. 
\end{problem}

\section*{Acknowledgements}

We are grateful to Ga\v{s}per Ko\v{s}mrlj for providing us with his software that computes game domination invariants. We acknowledge the financial support from the Slovenian Research Agency (research core funding No.\ P1-0297 and projects J1-9109, J1-1693, N1-0095, N1-0108).

\end{document}